\documentclass[
a4paper
,11pt
,oneside
]{amsart}
\usepackage{amssymb}
\usepackage{amsthm}  
\usepackage{amsmath} 
\usepackage{amscd}
\usepackage{url}  
\usepackage[all]{xypic}
\usepackage{color}

\newcommand{\omi}{\omega^{-1}}
\newcommand{\om}{\omega}
\newcommand{\ka}{\kappa}
\newcommand{\si}{\sigma}

\newcommand{\ba}{\mathcal{G}}  
\newcommand{\va}{\varphi}
\newcommand{\fg}{\mathfrak g}
\newcommand{\fp}{\mathfrak p}
\newcommand{\Ad}{{\rm Ad}}
\newcommand{\ked}{{\rm ker}(\partial^*)}
\newcommand{\imd}{{\rm im}(\partial^*)}

\newcommand{\Fl}{{\rm Fl}}
\newcommand{\id}{{\rm id}}
\newcommand{\na}{\nabla}
\newcommand{\U}{\Upsilon}
\newcommand{\gr}{{\rm gr}}
\newcommand{\lz}{[\![}
\newcommand{\pz}{]\!]}
\newcommand{\R}{\mathbb{R}}
\newcommand{\C}{\mathbb{C}}
\newcommand{\Le}{\mathcal{L}}

\newcommand{\At}{\mathcal{A}}
\newcommand{\Rho}{{\mbox{\sf P}}}

\newtheorem*{prop*}{Proposition}

\newtheorem*{thm*}{Theorem}

\newtheorem*{lemma*}{Lemma}
\newtheorem*{rem*}{Remark}
\theoremstyle{definition}
\newtheorem*{def*}{Definition}

\theoremstyle{remark}

\newtheorem*{exam*}{Example}

\newtheorem*{cor*}{Corollary}

\begin{document}
\title{Symmetries of parabolic contact structures}
\author{Lenka Zalabov\' a}
\address{Institute of Mathematics and Biomathematics, Faculty of Science, 
University of South Bohemia, Brani\v sovsk\' a 31, \v Cesk\' e Bud\v ejovice,
 370 05, Czech Republic}
\email{lzalabova@gmail.com}
\thanks{The author would like to mention very useful discussions with 
Andreas \v Cap during the work on this paper and the support of 
Eduard \v Cech Center for Algebra and Geometry for her participation in 
the conference The 30th Winter School Geometry and Physics, 2010.}
\begin{abstract}
We generalize the concept of locally symmetric spaces to parabolic 
contact structures. We show that symmetric normal parabolic 
contact structures are torsion--free and some types of them have to be locally 
flat. We prove that each symmetry given at a point with non--zero 
harmonic curvature is involutive. Finally we give restrictions on number 
of different symmetries which can exist at such a point.
\end{abstract}
\maketitle
Affine locally symmetric spaces are well known and studied objects 
from differential geometry. The classical definition says that a 
local symmetry at $x$ on a manifold $M$ with an affine connection $\na$ is 
a locally defined affine transformation such that $s_x(x)=x$ and 
$T_xs_x=-\id$ on $T_xM$. We can understand $\na$ as a geometric structure 
on $M$ such that the symmetry $s_x$ preserves this structure. There is 
a natural generalization of this concept: For a manifold with an 
arbitrary geometric structure, one can define a local symmetry as 
an automorphism of this geometric structure, which satisfies the two 
above conditions. Best known examples of this concept are 
Riemannian symmetric spaces, see \cite{KN2}, and projective symmetric 
spaces, see \cite{P}. This concept also generalizes nicely to 
geometric structures which can be described as $|1|$--graded 
parabolic geometries, see \cite{ja-srni06, ja-elsevier}.

In this article, we are interested in symmetries of contact manifolds
endowed with some additional structures which can be described 
as parabolic geometries, the so--called 
parabolic contact structures, see \cite{parabook}. 
Discussion of the Levi bracket implies that we cannot define a symmetry 
in the classical sense, see \cite{ja-srni06}. Motivated by the definition 
of a symmetry for Cauchy--Riemann structures from \cite{KZ}, we define 
a symmetry at $x$ as a morphism of the contact geometry such that 
$s_x(x)=x$ and $T_xs_x=-\id$ on the contact distribution at $x$. 
This definition works nicely for all parabolic contact structures. Then, 
to study symmetries on parabolic contact structures, we can use 
general techniques known from theory of parabolic geometries, 
see \cite{WS,parabook}. 

In this article, we discuss the curvature of symmetric parabolic 
contact geometries in detail. The theory of harmonic curvature for 
parabolic geometries allows us to prove that symmetric normal 
parabolic contact geometries must be torsion--free. Moreover, some types 
of them have to locally flat, if they are symmetric. 

Contrary to the classical case, symmetries of parabolic contact 
structures are not involutive in general and we use Weyl structures to 
study this question. More precisely, we show that only locally 
flat geometries can carry non--involutive symmetries at each point. 
Moreover, for each involutive symmetry on a parabolic contact geometry, 
there exists an admissible affine connection which is invariant with 
respect to the symmetry. Finally we show that in many cases, there can 
exist at most one symmetry at points with non--zero curvature.

\section{Parabolic contact structures}
We remind here basic definitions and summarize properties of 
parabolic geometries. We discuss parabolic contact structures in detail. 
We introduce  here Weyl structures which are our main tool to deal 
with parabolic contact geometries. We follow concepts and notation of 
\cite{parabook, WS} and the reader can find all details and proofs therein. 

\subsection{Contact structures and parabolic geometries} \label{contact}
Consider a manifold $M$ endowed with a distribution $H$ of $TM$ of 
corank one. Then $H \subset TM$ forms a filtration and on the graded 
bundle $\gr(TM)=H \oplus TM/H$, there is the 
\emph{Levi bracket} $\mathcal{L}:H \times H\rightarrow TM/H$ which is 
a bilinear bundle map induced by the Lie bracket of vector fields. The 
well known definition says that $H \subset TM$ forms a 
\emph{contact structure} on $M$ if the Levi bracket is non--degenerate 
at each point. The subbundle $H$ is then called \emph{contact distribution}. 

We will discuss here contact manifolds endowed with some additional structures which can be described as parabolic geometries. Let us remind that for a 
semisimple Lie group $G$ and its parabolic subgroup $P$, 
a~\emph{parabolic geometry} of type $(G,P)$ is a pair 
$(p:\ba\rightarrow M,\om)$ consisting of a principal $P$--bundle 
$\ba \rightarrow M$ and of a $1$--form $\om \in \Omega^1(\ba,\fg)$, 
called \emph{Cartan connection}, which is $P$--equivariant, 
reproduces generators of fundamental vector fields and induces a 
linear isomorphism $T_u\ba\cong\fg$ for each $u \in \ba$. The Lie algebra 
$\fg$ of $G$ is then equipped (up to the choice of Levi factor $\fg_0$ in 
$\fp$) with a grading of the form 
$\fg_{-k}\oplus \dots \oplus \fg_{0} \oplus \dots \oplus\fg_{k}$ 
such that the Lie algebra of $P$ is exactly 
$\fp = \fg_0 \oplus \dots \oplus \fg_k$. 
There is the usual notation such that 
$\fg_{-}:=\fg_{-k} \oplus \dots \oplus \fg_{-1}$, 
$\fp_+:=\fg_{1} \oplus \dots \oplus \fg_k$ and $P_+ \subset P$ is the 
subgroup corresponding to $\fp_+$. By $G_0$ we denote the subgroup in 
$P$, with Lie algebra $\fg_0$, consisting of all elements whose 
$\Ad$--action preserves the grading of $\fg$. Each element $g$ of $P$ can 
be uniquely written as $g_0\exp Z_1 \cdots \exp Z_k$ for suitable 
$g_0 \in G_0$ and $Z_i \in \fg_i$, thus $\exp Z_1 \cdots \exp Z_k \in P_+$. 
Let us remind that for each parabolic geometry, there is an element 
$E \in \fg_{0}$ with the property $[E,X]=iX$ for each $X \in \fg_{i}$, 
the so--called \emph{grading element}.
 To study contact structures, we have to focus on a special case of 
$|2|$--grading: A \emph{contact grading} of a simple Lie algebra $\fg$ is 
a grading 
$\fg_{-2} \oplus \fg_{-1} \oplus \fg_{0} \oplus \fg_{1} \oplus \fg_{2}$ 
such that $\fg_{-2}$ is one--dimensional and the Lie bracket 
$[- ,- ] : \fg_{-1} \times \fg_{-1} \rightarrow  \fg_{-2}$ is non--degenerate. 
Let us remark that for each contact grading,  the subspace 
$[\fg_{-2},\fg_{2}]$ coincides with the subspace generated by $E$.

It is well known that the Cartan connection $\om$ provides an 
identification $TM \simeq \ba \times_P \fg/\fp$. Suppose we have a 
parabolic geometry corresponding to a contact grading. Because each 
contact grading of $\fg$ induces $P$--invariant filtration of the 
form 
$\fg=\fg^{-2}\supset \fg^{-1} \supset \fg^0 \supset \fg^1 \supset \fg^{2}=\fg_2$, 
the subspace $\fg^{-1}/\fp \subset  \fg/\fp$ defines a subbundle 
$T^{-1}M := \ba \times_P \fg^{-1}/\fp$ of corank one in $TM$. There is 
the Levi bracket on $\gr(TM)=T^{-1}M \oplus TM/T^{-1}M$ and the geometry 
is called \emph{regular}, if the Levi bracket corresponds to the Lie 
bracket $[-,- ] : \fg_{-1} \times \fg_{-1} \rightarrow  \fg_{-2}$ under 
the above identification. Then for regular parabolic geometries 
corresponding to contact gradings, the underlying filtration 
$T^{-1}M \subset TM$ defines a contact structure on $M$, and each 
such geometry is called \emph{parabolic contact structure} or 
\emph{parabolic contact geometry}.
Moreover, define $\ba_0 := \ba/P_+$, which is a principal $G_0$--bundle 
over $M$. This is the reduction of the natural frame bundle of $\gr(TM)$ 
with respect to $\Ad: G_0 \rightarrow  Gl(\fg_{-1})$ and in this way, we 
get an additional geometric structure on $T^{-1}M$. 

Let us remind some more facts on parabolic contact structures that 
will be needed:
The $P$--bundle $G\rightarrow G/P$ together with the (left) 
Maurer--Cartan form $\om_G \in \Omega^1(G,\fg)$ forms a geometry which 
is called \emph{homogeneous model}. 
A \emph{morphism} between geometries of type $(G,P)$ from 
$(\ba \rightarrow M,\om)$ to $(\ba'\rightarrow M',\om')$ is a 
$P$--bundle morphism $\va:\ba \rightarrow \ba'$ such that $\va^*\om'=\om$.
 We will suppose that the maximal normal subgroup of $G$ which is contained 
in $P$ is trivial. With this assumption, there is one--to--one 
correspondence between morphisms of parabolic geometries and their 
base morphisms. Let us remind that such geometries are called \emph{effective}.

The curvature is described by $P$--equivariant mapping 
$\ka:\ba \rightarrow \wedge^2 (\fg/\fp)^* \otimes \fg$, 
the so--called \emph{curvature function}. The Maurer--Cartan equation 
implies that the curvature of the homogeneous model vanishes. Conversely, 
it can be proved that if the curvature of a geometry vanishes, then it 
is locally isomorphic to the homogeneous model of the same type. If $\ka$ 
has its values in a subbundle $\wedge^2 (\fg/\fp)^* \otimes \fp$, we call 
the geometry \emph{torsion--free}.
The regular geometry is called \emph{normal}, if the curvature 
satisfies $\partial^* \circ \ka=0$, where $\partial^*$ is the differential 
in the standard complex computing Lie algebra homology of $\fp_+$ 
with coefficients in $\fg$. Then we can define the 
\emph{harmonic curvature} $\ka_H$ which is the composition of the 
curvature function with  the projection $\ked \rightarrow \ked/\imd$. 
There is the following general statement, see \cite{parabook}:
\begin{thm*} \label{harmonic}
On a regular normal parabolic geometry, the curvature $\ka$ vanishes 
if and only if the harmonic curvature $\ka_H$ vanishes.
\end{thm*}
It can be proved that $\ked/\imd$ is a $G_0$--submodule of 
$\wedge^{2}\fg_{-}^* \otimes \fg$ and decomposes into the direct sum 
of components each of which is contained in some 
$\fg_{-i}^* \wedge \fg_{-j}^* \otimes \fg_k$. 
According to this decomposition, $\ka_H$ decomposes into the sum of 
components of homogeneity $\ell=i+j+k$ that we denote by $\ka^{(\ell)}$.
One can use the Kostant's version of the Bott--Borel--Weil theorem to find 
all components of $\ka_H$.
 For parabolic contact geometries, it turns out that there can exist 
only the following three types of components:
\begin{itemize}
\item $\ka^{(1)}$ valued in $\fg_{-1}^* \wedge \fg^*_{-1} \otimes \fg_{-1}$, 
\item  $\ka^{(2)}$ valued in $\fg_{-1}^* \wedge \fg^*_{-1} \otimes \fg_{0}$ and
\item $\ka^{(4)}$ valued in $\fg_{-2}^* \wedge \fg^*_{-1} \otimes \fg_{1}$. 
\end{itemize}
See Appendix for more detailed description of all contact gradings 
and corresponding geometries with their components of harmonic curvature.

\subsection{Adjoint tractor bundles and Weyl structures} \label{WeylStructures}
Let us introduce here briefly concept of adjoint tractor bundles which 
allows us to write formulas and make computations in a more convenient form. 
 The \emph{adjoint tractor bundle} is the natural bundle 
$\At M := \ba \times_P \fg$ corresponding to the restriction of 
$\Ad$--action of $G$ on $\fg$.  
For each parabolic contact geometry,
the filtration of $\fg$ induces a filtration 
$\At M=\At^{-2}M\supset \At^{-1}M \supset \At^0M \supset \At^1M \supset \At^{2}M$ 
such that $\At^{i} M=\ba \times_P \fg^i$, and there is the associated 
graded bundle 
$\gr(\At M)=\At_{-2}M \oplus \At_{-1}M \oplus \At_0M \oplus \At_1M \oplus \At_{2}M$, 
where $\At_iM=\At^iM/\At^{i+1}M$ equals to $\ba_0 \times_{G_0} \fg_i$. 
 Clearly, $TM \simeq \At M/\At^0 M$ and $T^*M\simeq\At^1 M$. On the 
graded bundle $\gr(\At M)$, there is the \emph{algebraic bracket} 
$ \{- ,-  \}: \At_{i}M \times \At_{j}M \rightarrow \At_{i+j}M$ defined 
by means of the Lie bracket on $\fg$. 
Clearly, its part $\At_{-1}M \times \At_{-1}M \rightarrow \At_{-2}M$ 
on $\gr(TM)=\At_{-2}M \oplus \At_{-1}M$ coincides with the Levi bracket 
thanks to the regularity. Since each fiber of $\gr(\At M)$  is isomorphic 
to $\fg$, the grading element defines a unique element 
$E(x) \in \gr_0(\At_xM)$ such that $\{E(x),- \}$ is a multiplication by $i$ 
on $\gr_i(\At_xM)$. In fact, these elements form a section $E$ of 
$\gr_0(\At M) $ which is called the \emph{grading section}. 
Let us also remark that we simultaneously get an action $\bullet$ of 
$\gr(\At M)$ on arbitrary tensor products of $\gr(\At M)$ which is given 
using the tensoriality of the algebraic bracket. In particular, the 
grading section acts on each homogeneous component of the tensor 
as multiplication by its homogeneity. 

Now we should remind basic facts on Weyl structures. For any parabolic 
contact geometry $(\ba \rightarrow M, \om)$ with the underlying 
$G_0$--bundle $p_0: \ba_0 \rightarrow M$, a \emph{Weyl structure} is a 
global smooth $G_0$--equivariant 
section $\si:\ba_0\rightarrow\ba$ of the canonical 
projection~$\pi:\ba \rightarrow \ba_0$. Weyl structures always exist and 
for any two Weyl structures $\si$ and $\hat \si$, there are 
$G_0$--equivariant functions 
$\U_1:\ba_0\rightarrow \fg_1$ and $\U_2:\ba_0\rightarrow \fg_2$ such that 
\begin{align*} 
\hat \si(u_0)= \si(u_0) \exp \U_1(u_0) \exp \U_2(u_0) 
\end{align*}
for all $u_0 \in \ba_0$. Clearly, $\U_i \in \Gamma(\At_iM)$ 
and $\U:=(\U_1,\U_2)$ is a smooth section of $\gr(T^*M)$.  Moreover, 
the Campbell--Baker--Hausdorff formula implies 
$\exp \U=\exp \U_1 \exp \U_2=\exp(\U_1+\U_2)= \exp \U_2 \exp \U_1$, 
see \cite{KMS}. 

For each Weyl structure $\si$, we can form the pullback
 $\si^*\om \in \Omega^1(\ba_0,\fg)$. This decomposes as 
$\si^*\om=\si^*\om_{-}+\si^*\om_{0}+\si^*\om_{+}$ and the part
$\si^*\om_{-} \in \Omega^1(\ba_0,\fg_{-})$ is called the \emph{soldering form}. 
Each Weyl structure $\si$ induces by means of its soldering form 
an isomorphism $TM \simeq \gr(TM)$ which we write as 
$\xi \mapsto (\xi_{-2},\xi_{-1})$. If $\si \exp \U_1 \exp \U_2$ is 
another Weyl structure, the isomorphism changes as 
$\xi \mapsto (\xi_{-2},\xi_{-1}-\{\U_{1},\xi_{-2} \} )$. 
In particular, $\si$ and $\si \exp \U_2$ induce the same isomorphism for 
an arbitrary $\U_2:\ba_0 \rightarrow \fg_2$.

The part $\si^*\om_0 \in \Omega^1(\ba_0,\fg_0)$ defines a principal 
connection on $p_0:\ba_0 \rightarrow M$ which we call the 
\emph{Weyl connection}. 
This connection induces connections on all associated bundles. 
In particular, for each $\si$ we get a preferred connection on 
$\gr(TM)=\ba_0 \times_{G_0} \fg_{-}$ and via the above isomorphism, we get 
a preferred connection on the tangent bundle, cotangent bundle and 
their tensor products. We call each such connection Weyl connection, too. 
For a Weyl structure $\si$, we denote the corresponding connection 
by $\na^\si$. 
For $\si$ and $\hat \si=\si \exp \U_1 \exp \U_2$ we have 
\begin{align} \label{zmena-konexe}
\na^{\hat \si}_\xi s = \na^{\si}_\xi s + \big( {1 \over 2} \{\U_1, \{\U_1,\xi_{-2} \}\} - \{\U_2, \xi_{-2} \} - \{\U_1 ,\xi_{-1} \} \big) \bullet s,
\end{align}
where $\xi \in \frak X(M)$ and $s$ is a section of an appropriate bundle.
The positive part $\si^*\om_+ \in \Omega^1(\ba_0,\fp_+)$ is called 
\emph{Rho--tensor} and is denoted by $\Rho^\si$. We will not need 
it explicitly, see \cite{parabook, WS} for details.

Let us finally remind the so--called normal Weyl structures. 
A \emph{normal Weyl structure} at $u$ is the only $G_0$--equivariant 
section $\si_u:\ba_0 \rightarrow  \ba$ satisfying
$
\si_u \circ \pi \circ {\rm Fl}_1^{\omi(X)}(u) = {\rm Fl}_1^{\omi(X)}(u)
,$ 
where by $\Fl^{\omi(X)}_t(u)$ we denote flows of constant vector fields 
$\omi(X) \in \frak X(\ba)$.
Each normal Weyl structure $\si_u$ is defined locally over some 
neighborhood of $p(u)$ and depends only on the $G_0$--orbit of 
$u \in \ba$, see \cite{parabook}.

\section{Basic facts on symmetries}
We formulate here the definition of a symmetry on a parabolic contact 
geometry and  describe its basic properties. We study the action of 
symmetries on Weyl structures and we describe some interesting subclasses 
of them. We focus here on the question of involutivity of our symmetries.

\subsection{Definitions} \label{definitions}
Let $(\ba \rightarrow M, \om)$ be a parabolic contact structure.
A (\emph{local}) \emph{symmetry with the center at $x \in M$} is a 
(locally defined) diffeomorphism $s_x$ on $M$ such that:
\begin{enumerate} 
\item $s_x(x)=x$, 
\item $T_x s_x=-\id$ on $T_x^{-1}M$,
\item $s_x$ is a base morphism of some (locally defined) automorphism $\va$
of the parabolic contact geometry.
\end{enumerate}
The geometry is called (\emph{locally}) \emph{symmetric} 
if there is a (local) symmetry at each point $x \in M$.

Clearly, each symmetry is a local symmetry.
In this article, we discuss local symmetries and local properties of locally 
symmetric geometries and we will shortly say `symmetry at $x$' and `symmetric'
instead of `local symmetry at $x$' and `locally symmetric', respectively. 
Global symmetries and their systems we will discuss elsewhere. 
Moreover, we will also call the automorphism $\va$ of $\ba$ and its underlying automorphism $\va_0$ of $\ba_0$  a `symmetry at $x$'.

\subsection{Basic properties of symmetries} \label{basic}
Let $s_x$ be a symmetry on a parabolic contact geometry and let $\va$ be 
as above. Since each symmetry $s_x$ preserves $x$, the (uniquely 
given) automorphism $\va$ has to preserve the fiber over $x$. Then for 
each frame $u \in p^{-1}(x)$ we have 
$\va(u)=ug_0\exp Z=ug_0\exp Z_1 \exp Z_2$ for suitable 
$g_0 \in G_0$, $Z_1 \in \fg_1$ and $Z_2 \in \fg_2$, where $Z=Z_1+Z_2$. 
Let us describe the element $g_0\exp Z_1 \exp Z_2$ in detail: 

For each $\xi(x)=\lz u,X \pz \in T^{-1}_xM$, i.e. for each 
$X \in \fg^{-1}/\fp$, we have
\begin{align}
\begin{split} \label{prvek}
T_xs_x.\xi(x)&=\lz \va(u),X\pz=\lz ug_0\exp Z_1\exp Z_2 ,X \pz= \\ &=
\lz u, \underline{\Ad}_{\exp (-Z_2)}\underline{\Ad}_{\exp (-Z_1)} \Ad_{g_0^{-1}}X \pz
\end{split}
\end{align}
and simultaneously $T_xs_x.\xi(x)=-\xi(x)=\lz u, -X \pz$. All together, 
the element $g_0\exp Z_1 \exp Z_2$ has to induce $-\id$ on $\fg^{-1}/\fp$ 
by the $\underline \Ad$--action. Moreover, $\exp Z_1 \exp Z_2$ acts 
trivially on $\fg^{-1}/\fp$. Indeed, there is the formula 
\begin{align} \label{Ad-exp}
\Ad_{\exp Z}X=\sum_{j=0}^\infty{1  \over j!} {\rm ad}^j_{Z}X =X +[Z,X]+{1 \over 2} [Z,[Z,X]]+ \cdots
\end{align}
for all $X \in \fg_{-}$ and $Z \in \fp_+$ and if $X \in \fg_{-1}$, 
all brackets on the right hand side belong to $\fp$.  In fact, along the 
fiber over $x$, the $P_+$--parts of the above elements are determined by 
$\va$ and can be arbitrary in general, one only has to impose 
the compatibility of $\va$ with the right action of $P$. 
Then the element $g_0$ has to cause the sign change on 
$\fg^{-1}/\fp \simeq \fg_{-1}$. Since our geometries are effective, there 
can exist at most one element $g_0 \in G_0$ which gives $-\id$ on 
$\fg_{-1}$ and it has to be the same element along the fiber. In 
particular, the underlying morphism $\va_0$ is of the form 
$\va_0(u_0)=u_0g_0$ for each $u_0 \in p_0^{-1}(x)$. Clearly, the element 
$g_0$ has to induce identity on $\fg_{-2}=[\fg_{-1},\fg_{-1}]$. 

One of basic properties of classical symmetries is their involutivity 
and there is a natural question on involutivity of our symmetries. 
Thus let us focus on $s_x^2:=s_x \circ s_x$. Let us first point out that  
$\va \circ \va = \id_\ba$ if and only if $s_x \circ s_x = \id_M$, 
which follows directly from effectivity. Clearly, $\va_0$ is then 
involutive, too.
Thus it suffice to study the morphism $\va^2$. In the above 
notation, $\va^2(u)=ug_0\exp Z_1 \exp Z_2g_0\exp Z_1 \exp Z_2$ holds and
using the known fact $\exp X g_0=g_0  \exp(\Ad_{g_0}X)$, we can rewrite 
this as $ug_0^2 \exp (-Z_1) \exp Z_2\exp Z_1 \exp Z_2.$ Moreover, $g_0^2$ 
acts as $\id$ on $\fg_{-1}$ and thus on $\fg$, because $\fg_{-1}$ 
generates $\fg_{-}$, $\fp =\fg_{-}^*$ and 
$\fg_{0} \subset \fg_{-}^* \otimes \fg_{-}$. 
Thus it lies in the kernel of the $\Ad$--action which coincides with 
the maximal normal subgroup of $G$ which is contained in $P$. 
Effectivity then gives $g_0^2=e$.
All together, we have got $\va^2(u)=u \exp 2Z_2.$ In another frame $uh$ 
for $h \in P$ we then have $\va^2(uh)= uh  \exp 2\Ad_{h}Z_2$.
Thus we can simply view $Z_2$ as $P$--equivariant function 
$Z_2: p^{-1}(x) \rightarrow \fg_2$. In fact, we have proved the 
following statement:
\begin{lemma*}
On a parabolic contact geometry, each symmetry $s_x$ at $x$ defines 
uniquely a covector $Z_2 \in T^*_xM$
through the equation 
\begin{align} \label{vau}
\va^2(u)=u \exp 2Z_2(u)
\end{align} 
along the fiber over $x$. The symmetry $s_x$ is involutive if and only if 
the covector $Z_2$ equals to zero.
\end{lemma*}
Thanks to the above observations, it is easy to describe the differential 
of each symmetry at its center.
\begin{prop*}
For each symmetry $s_x$ at $x$ on a parabolic contact geometry, the 
mapping $T_xs_x:T_xM \rightarrow T_xM$ is involutive,
thus $T_xM$ decomposes into two eigenspaces with eigenvalues $-1$ and $1$. 
The eigenspace corresponding to the eigenvalue $-1$ has to coincide 
with $T^{-1}_xM$, the contact distribution, and there exists 
one--dimensional eigenspace corresponding to the eigenvalue $1$. 
\end{prop*}
\begin{proof}
For each $\xi(x)=\lz u, X \pz$ from $T_xM$ we have 
\begin{align*}
T_xs_x^2.\xi(x)&= \lz \va^2(u), X\pz = \lz u \exp 2Z_2(u),X \pz =
\lz u, \underline \Ad_{\exp (-2Z_2(u))}X \pz \\&=\lz u, X \pz=\xi(x) 
,
\end{align*}
which follows directly from formulas $(\ref{vau})$ and $(\ref{Ad-exp})$. 
The rest follows immediately from the definition of the symmetry.
\end{proof}

\subsection{Action of symmetries on Weyl structures} \label{action}
Let us now discuss relations of various Weyl structures to the symmetry $s_x$. 
For each Weyl structure $\si$ we can write
\begin{align} \label{obecne}
\va(\si(u_0))=\si(\va_0(u_0)) \exp \U_1(u_0) \exp \U_2(u_0)
\end{align}
for each $u_0 \in \ba_0$ and for suitable functions 
$\U_1: \ba_0\rightarrow \fg_1$ and $\U_2: \ba_0\rightarrow \fg_2$ which 
are generally determined by $\va$ and $\si$. 
With the notation from the last section, we have $\U_2(u_0)=Z_2(\si(u_0))$ 
in the fiber over $x$. The Lemma \ref{basic} then shows that $\U_2$ does 
not depend on the choice of a Weyl structure $\si$ and coincides for all 
Weyl structures at $x$. Clearly, $\U_2$ vanishes at $x$ if and only if 
$s_x$ is involutive.
The function $\U_1$ depends on the choice of a Weyl structure $\si$ at $x$ 
and with the above notation, $\U_1(u_0)=Z_1$ for $\si(u_0)=u$.

Let us now focus on the role of $\U_1$ for the isomorphism 
$TM \simeq \gr(TM)$ given by the Weyl structure $\si$: For a tangent 
vector $\xi(x)=\lz \si(u_0),X_{-2}+X_{-1} \pz$ where $X_i \in \fg_i$ we have 
\begin{align*}
T_xs_x.\lz \si(u_0),X_{-2}+X_{-1} \pz&=
\lz \si(\va_0(u_0)) \exp \U_1(u_0) \exp \U_2(u_0),X_{-2}+X_{-1}\pz \\&=
\lz \si(u_0), X_{-2}-X_{-1} - [\U_1(u_0), X_{-2}]\pz,
\end{align*}
which follows from  the fact that $\va_0(u_0)=u_0g_0$ for $g_0$ giving 
$-\id$ on $\fg_{-1}$ and from formulas $(\ref{obecne})$ and $(\ref{Ad-exp})$.
In particular, the isomorphism $TM \simeq \gr(TM)$ given by a Weyl 
structure $\si$ reflects the decomposition of $T_xM$ into 
$\pm 1$--eigenspaces for $T_xs_x$ if and only if the Weyl structure 
$\si$ satisfies $\U_1(u_0)=0$ for each $u_0$ from the fiber over $x$. 
\begin{lemma*}
On a parabolic contact geometry with a symmetry $s_x$ at $x$, there are 
Weyl structures $\hat \si$ such that 
$
\va(\hat \si(u_0))=\hat \si(\va_0(u_0)) \exp \hat \U_1(u_0) \exp  \U_2(u_0)
$
holds for suitable $\hat \U_1$ such that
$\hat \U_1(u_0)=0$ for each $u_0$ from the fiber over $x$.
\end{lemma*}
\begin{proof}
Consider an arbitrary Weyl structure $\si$ and let $(\U_1, \U_2)$ 
be determined by $\si$ as above. Let us verify that the Weyl structure
\begin{align*} 
\hat \si(u_0)=\si(u_0) \exp (-{1 \over 2}\U_1(u_0)) 
\end{align*} 
satisfies the condition: The formula $(\ref{obecne})$ and 
the Campbell--Baker--Hausdorff formula allow us to write
\begin{align*}
\va(\hat \si(u_0)) &= \va(\si(u_0)) \exp (-{1 \over 2}\U_1(u_0))
\\&= \si(\va_0(u_0))  \exp \U_1(u_0)  \exp \U_2(u_0) 
\exp( -{1 \over 2}\U_1(u_0))
\\&=\si(\va_0(u_0))  \exp{1 \over 2}\U_1(u_0)  \exp \U_2(u_0). 
\end{align*}
Equivariancy of $\U_1$ gives $\U_1(\va_0(u_0)) = \U_1(u_0g_0)=-\U_1(u_0)$ 
in the fiber over $x$ for $g_0$ giving $-\id$ on $\fg_{-1}$ and we can 
rewrite the above expression as  
$$
\si(\va_0(u_0))  \exp (-{1 \over 2}\U_1(\va_0(u_0)))  \exp\U_2(u_0)=
\hat \si(\va_0(u_0))  \exp\U_2(u_0)
$$ 
in the fiber over $x$. Thus $\hat \si$ is the required Weyl structure. 
\end{proof}
Let us call each Weyl structure $\hat \si$ satisfying the condition in the 
Lemma an \emph{almost $s_x$--invariant Weyl structure} at $x$. All almost 
$s_x$--invariant Weyl structures form a family of Weyl structures which 
is parametrized over $\fg_2$ at $x$. Really, all Weyl structures inducing 
the same isomorphism $T_xM\simeq \gr(T_xM)$ as $\hat \si$ are of the 
form $\hat \si \exp F_1(u_0)\exp F_2(u_0)$ for arbitrary functions 
$F_2: \ba_0 \rightarrow \fg_2$ and $F_1:\ba_0 \rightarrow \fg_1$ 
where $F_1(u_0)=0$ in the fiber over $x$, see \ref{WeylStructures}. 

Let us finally describe the involutivity of our symmetries in the language 
of  Weyl structures.
\begin{prop*}
On a parabolic contact geometry with a symmetry $s_x$ at $x$,
the following facts are equivalent:
\begin{enumerate}
\item[$(a)$] the symmetry $s_x$ is involutive,
\item[$(b)$] there exists a Weyl structure $\si$ such that 
$
\va(\hat \si(u_0))=\si(\hat \va_0(u_0))
$
holds in the fiber over $x$,  
\item[$(c)$] there exists a Weyl structure $\si_u$ such that 
$
\va( \si_u(u_0))=\si_u(\va_0(u_0))
$
holds over some neighborhood of $x$.
\end{enumerate}
\end{prop*}
\begin{proof}
$(a) \Rightarrow (b)$
Let $\hat \si$ be an arbitrary almost $s_x$--invariant Weyl structure. 
The Lemma \ref{basic} says that the involutivity 
implies vanishing of $\U_2$ in the fiber over $x$.  
Thus if $s_x$ is involutive, the almost $s_x$--invariant Weyl structure 
$\hat \si$ has to satisfy $(b)$.
\smallskip \\
$(b) \Rightarrow (c)$
Let $\hat \si$ be an arbitrary Weyl structure satisfying 
$\va(\hat \si(u_0))=\hat \si(\va_0(u_0))$
in the fiber over $x$. Consider the normal Weyl structure $\si_u$ such that 
$\si_u(u_0) = \hat \si(u_0)$ for $p_0(u_0)=x$. The condition of the normality 
prescribes $\si_u$ uniquely on a normal neighborhood of $x \in M$, 
see \ref{WeylStructures} for definition. But then, because 
$\va(\si_u(u_0))= \si_u(\va_0(u_0))$  
holds in the fiber over $x$, it has to hold over some normal neighborhood of $x$
and $\si_u$ satisfies $(c)$.
\smallskip \\
$(c)\Rightarrow (a)$
Consider an arbitrary Weyl structure $\si$ satisfying $(c)$. This can 
be equivalently written as $\va^{-1}(\si(\va_0(u_0)))=\si(u_0)$ which 
means that the corresponding Weyl connection is invariant with respect to 
$s_x$. Since the isomorphism $T_xM \simeq \gr(T_xM)$ reflects 
the decomposition of $T_xM$ into $\pm 1$--eigenspaces, we can describe 
$s_x$ on a neighborhood of $x$ nicely via geodesics of the 
invariant connection. Indeed, each vector 
$(\xi_{-2}(x),\xi_{-1}(x)) \in T_xM$ determines uniquely a geodesic at 
$x$, and the symmetry $s_x$ maps it on a geodesic at $x$, which is 
uniquely determined by a vector $(\xi_{-2}(x),-\xi_{-1}(x))$. This 
describes $s_x$ on a neighborhood of $x$ and one can see directly that 
it has to be involutive.
\end{proof} 
Let us call each Weyl structure satisfying the condition $(b)$ of 
the Proposition  \emph{$s_x$--invariant Weyl structure at $x$} and each 
Weyl structure satisfying the condition $(c)$ of the Proposition 
\emph{$s_x$--invariant Weyl structure on a neighborhood of $x$}. 

\section{Symmetries of homogeneous models} \label{homogeneous}
In this section, we focus on homogeneous models, which are simplest 
examples of parabolic contact geometries. We describe explicitly 
their symmetries and we give some concrete examples of homogeneous 
symmetric geometries. 

\subsection{Description of symmetries}
Let $(G\rightarrow G/P,\om_G)$ be a homogeneous model of a parabolic 
contact geometry of type $(G,P)$. It is well known that all automorphisms 
of the homogeneous model are 
just left multiplications by elements of $G$ and an analog of the 
Liouville theorem 
states that any local automorphism can be uniquely extended to a global one, 
see \cite{S,parabook}.  Thus if the homogeneous model is locally symmetric, 
then it is symmetric. Moreover, because $G$ acts transitively on $G/P$,
it suffices to find a symmetry at the origin to decide whether the 
homogeneous model is symmetric. 
\begin{prop*}
All symmetries of the homogeneous model of a parabolic contact geometry 
centered at the origin $o=eP$ are given by left multiplications  
by elements $g_0 \exp Z_1 \exp Z_2 \in P$, where $Z_1\in\fg_1$ 
and $Z_2\in\fg_2$ 
are arbitrary and $g_0 \in G_0$ is such that $\Ad_{g_0}=-{\rm id}$ 
on $\fg_{-1}$. 
In particular, if there is one symmetry at a point, then there is an infinite 
amount of them.
\end{prop*} 
\begin{proof}
For homogeneous models,
$
T^{-1}(G/P) = G \times_P \fg^{-1}/\fp.
$
Then we can write each tangent vector $\xi(o) \in T^{-1}_o(G/P)$ as
$
\xi(o)=\lz e, X \pz
$
for suitable $X \in \fg^{-1}/\fp$. Since automorphisms of the homogeneous model 
are left multiplications $\lambda_g$ by elements  $g \in G$, all symmetries 
at the origin are exactly left multiplications $\lambda_g$ satisfying
$\lambda_g(o)=o$ and $ T_o\lambda_g.\xi(o)=-\xi(o)$ for all contact 
vectors $\xi(o)$.
The first condition is equivalent to the fact that $g \in P$. Then $g$ can 
be written as
$g=g_0 \exp Z_1 \exp Z_2$ and the second condition means that
$$
T_o\lambda_{g_0 \exp Z_1 \exp Z_2}.\lz e,X\pz= 
\lz g_0 \exp Z_1 \exp Z_2,X\pz=
\lz e, \underline{\Ad}_{g_0 \exp Z_1 \exp Z_2}^{-1} X \pz
$$
and 
$
-\xi(o)=\lz e, -X \pz
$
coincide for each $X \in \fg^{-1}/\fp$. Thus we look for elements 
$g \in P$ such that 
$ \underline{\Ad}_{\exp(-Z_2)} \underline{\Ad}_{\exp(-Z_1)} \Ad_{g_0^{-1}}X=-X$
for all $X \in \fg^{-1}/\fp$ and the rest follows immediately 
from observations in Section \ref{basic}.
\end{proof}
Let us finally discuss involutivity of these symmetries. The symmetry 
$g_0 \exp Z_1 \exp Z_2$ is involutive if and only if the element 
$(g_0 \exp Z_1 \exp Z_2)^2$ induces identity on $G/P$ and effectivity 
says that it has to be equal to $e$. We have
$$
g_0 \exp Z_1 \exp Z_2 g_0 \exp Z_1 \exp Z_2=e\exp 2 Z_2
.
$$
Thus involutive symmetries at the origin are left multiplications by 
elements $g_0 \exp Z_1$ where $g_0$ and $Z_1$ are as above and $Z_2$ has to 
be equal to zero. If $Z_2$ is non--zero,
then the symmetry is not involutive. In particular, there 
exist non--involutive symmetries on homogeneous models.

Clearly, if $g$ induces (involutive) symmetry at the origin $o = eP$, then $hgh^{-1}$ induces (involutive) symmetry at the point $hP$.

\subsection{Examples} \label{examples}
Let us introduce here some examples of parabolic contact structures  
and discuss symmetries on their homogeneous models, see \cite{parabook} 
for detailed description.

\subsubsection*{Lagrangean contact structures.} 
Let us start with $\fg=\frak{sl}(n+2,\R)$, the split real form 
of $\frak{sl}(n+2,\C)$, for $n \geq 1$. This admits a contact grading which 
is given by the following decomposition into blocks of sizes $1$, $n$ and $1$:
$$
\left( \begin{smallmatrix}
\fg_{0} & \fg_{1}^L & \fg_2 \\ \fg_{-1}^L & \fg_0 & \fg_{1}^R \\ \fg_{-2} & \fg_{-1}^R & \fg_0 \end{smallmatrix} \right).
$$
The splittings $\fg_{\pm 1}=\fg_{\pm 1}^L \oplus \fg_{\pm 1}^R$ are 
$\fg_0$--invariant and $\fg_{-1}^L$ and $\fg_{-1}^R$ are isotropic for 
$[- ,- ]:\fg_{-1} \times \fg_{-1} \rightarrow \fg_{-2}$. 
Let us choose $G=PGL(n+2,\R)$, the quotient of $GL(n+2,\R)$ by its 
center. Then $P$ consists of classes of block upper triangular matrices 
and $G_0$ of block diagonal matrices. In particular, $G_0$ coincides by 
means of the $\Ad$--action with the group of all automorphisms of graded 
Lie algebra $\fg_{-2} \oplus \fg_{-1}$ which in addition preserve 
the decomposition $\fg_{-1}=\fg_{-1}^L \oplus \fg_{-1}^R$. 
Thus for a parabolic contact geometry of type $(G,P)$, the underlying 
geometry consists of a contact distribution together with a 
fixed decomposition of the form $T^{-1}M=L \oplus R$ into two subbundles 
(of rank $n$) each of which is isotropic with respect to $\Le$. 
These geometries are known as Lagrangean contact structures.  
The homogeneous model is the flag manifold of lines in hyperplanes 
in $\R^{n+2}$.

Let us now discuss symmetries at the origin of the homogeneous model.
We look for an element $g_0 \in G_0$ such that $\Ad_{g_0}X=-X$ 
for each $X \in \fg_{-1}$.
Elementary matrix computation shows that there is a solution which 
is represented by the matrix of the form
$$
g_0=\Big( \begin{smallmatrix}
-1 & 0 & 0 \\0 & E & 0 \\ 0 & 0 & -1
\end{smallmatrix} \Big)
,$$
where $E$ is the identity matrix, 
and thus the homogeneous model is symmetric. All symmetries at the origin 
are represented by matrices of the form 
$$
\left( \begin{smallmatrix}
-1 & -V & \gamma  \\ 0 & E & W \\ 0 & 0 & -1
\end{smallmatrix} \right),
$$ 
where $V^*,W \in \R^{n}$ and $\gamma \in \R$ are arbitrary, and the 
involutive ones have to satisfy $\gamma=-{1 \over 2}VW$.

\subsubsection*{Non--degenerate partially integrable almost CR--structures 
of hypersurface type} 
Consider the real form $\fg=\frak{su}(p+1,q+1)$ of $\frak{sl}(n+2,\C)$ 
for $p+q=n\geq 1$. For a suitable choice of Hermitian product, this admits 
a contact grading of the same block form as the Lagrangean case. 
Denoting $\mathbb{I}$ the diagonal $n\times n$--matrix with the first 
$p$ entries equal to $1$ and the remaining $q$ entries equal to $-1$, 
we  write the elements explicitly~as 
$$
\left( \begin{smallmatrix}
a & Z & iz \\ X & A & -\mathbb{I}Z^* \\ ix & -X^*\mathbb{I} & -\bar a
\end{smallmatrix} \right)
,$$
where $x,z \in \R$, $a \in \C$, $X,Z^{T} \in \C^{n}$, $A \in \frak{u}(n)$ 
and $a +{\rm tr}A-\bar a=0$. The bracket
$\fg_{-1} \times \fg_{-1} \rightarrow \fg_{-2}$ is given by 
$[X, Y] = Y^*\mathbb{I}X-X^*\mathbb{I}Y$, which is twice the imaginary part 
of the standard Hermitian product of signature $(p,q)$.
Choose $G=PSU(p+1,q+1)$, the quotient of $SU(p+1,q+1)$ by its center. 
Then $P$ consists of classes of block upper triangular matrices and 
elements of $G_0$ are represented by block diagonal matrices. For 
parabolic contact geometries of type $(G,P)$, the underlying geometry 
consists of a contact distribution $T^{-1}M$ together with a complex 
structure $J$ such that $\Le(J\xi,J\eta)=\Le(\xi,\eta)$ for all 
$\xi,\eta \in \Gamma( T^{-1}M)$. 
Such geometries are known as partially integrable almost CR--structures 
of hypersurface type. The homogeneous model is the projectivized null cone 
of a Hermitian form of signature $(p+1,q+1)$, which is a real hypersurface 
in $\C P^{n+1}$.

Let us now discuss briefly symmetries of the homogeneous model at the origin.
We look for an element $g_0 \in G_0$ such that
$\Ad_{g_0}X=-X$ for all $X \in \fg_{-1}$. Elementary computation shows that 
the solution exists and is given by the same matrix $g_0$ as in the 
Lagrangean case. Then all symmetries at the origin are represented by 
matrices of the form 
$$
\left( \begin{smallmatrix}
-1 & -Z & iz \\ 0 & E & -\mathbb{I}Z^* \\ 0  & 0 & -1
\end{smallmatrix} \right)
$$
where $Z \in \mathbb{C}^{n*}$ and $z \in \R$ are arbitrary, and the 
involutive ones have to satisfy $iz=Z\mathbb{I}Z^*$.

\section{Curvature restrictions}
In this section, we discuss restrictions on the curvature of a 
parabolic contact geometry, which are caused by the existence of a 
symmetry. We study the torsion of symmetric parabolic contact geometries 
in detail. We show that there are relations between the curvature of 
a symmetric geometry and involutivity of its symmetries.

\subsection{Torsion restrictions}
Let us work with normal parabolic contact geometries here. In fact, 
the  normality assumption is only some technical restriction which plays 
no role, if we understand symmetries as morphisms of the underlying 
geometry, and this clearly is the most reasonable point of view. For 
such underlying geometry, there are various non--isomorphic 
parabolic geometries inducing this underlying structure and it can be 
proved that the normal one always exists, see \cite{parabook}. Assuming 
normal geometry, we can discuss its components of harmonic curvature, 
which are easily computable and provide an information on the whole 
curvature of the parabolic geometry, see \ref{contact} and Appendix. Let 
us start with the harmonic torsion $\ka^{(1)}$, which has its values 
in $\fg_{-1}^* \wedge \fg_{-1}^*\otimes \fg_{-1}$.
\begin{lemma*}
If there is a symmetry $s_x$ at $x$ on a normal parabolic contact geometry, 
then $\ka^{(1)}$ vanishes at $x$.
\end{lemma*} 
\begin{proof}
Let $\va$ be as usual. For $u \in p^{-1}(x)$ we have $\va(u)=ug_0 \exp Z$ 
for suitable $g_0 \in G_0$ and $Z \in \fp_+$, see \ref{basic}. Then for 
each $X,Y \in \fg_{-1}$ we get 
$$
\ka^{(1)}(\va(u))(X,Y)=\ka^{(1)}(u g_0 \exp Z)(X,Y)=
\exp(-Z) g_0^{-1} \cdot \ka^{(1)}(X,Y)$$ 
where $\cdot$ denotes the induced $\Ad$--action on 
$\fg_{-1}^* \wedge \fg_{-1}^*\otimes \fg_{-1}$.
The action of $g_0^{-1}$ is of the form
$$g_0^{-1} \cdot \ka^{(1)}(X,Y)=
\Ad_{g_0^{-1}}(\ka^{(1)}(u)(\Ad_{g_0}X,\Ad_{g_0}Y))=
$$
$$
-\ka^{(1)}(u)(-X,-Y)=-\ka^{(1)}(u)(X,Y)
$$
since the element $g_0$ acts as $-\id$ on $\fg_{-1}$, and the action 
of $\exp(-Z)$ is trivial.
Because automorphisms preserve curvature, $-\ka^{(1)}(u)(X,Y)$ has to be 
equal to $\ka^{(1)}(u)(X,Y)$ in the fiber over $x$ and then it has 
to vanish at $x$. 
\end{proof}
\begin{prop*}
Each normal symmetric parabolic contact geometry is torsion--free. 
Moreover, normal symmetric
\begin{itemize}
\item  Lie contact structures,
\item  parabolic contact geometries corresponding to exotic Lie algebras
\end{itemize}
have to be locally isomorphic to the homogeneous models of the same type.
\end{prop*}
\begin{proof}
Thanks to the regularity, the curvature satisfies 
$\ka(u)(\fg^i,\fg^j) \subset \fg^{i+j+\ell}$ for all $u \in \ba$ and for 
some $\ell \geq 1$ in general. Moreover, it can be proved that the 
component of degree $\ell$ mapping $\fg_i \times \fg_j$ to 
$\fg_{i+j+\ell}$ corresponds to the component of $\ka_H(u)$ of degree 
$\ell$, see \cite{parabook}. The above Lemma shows that $\ell \geq 2$ 
for symmetric geometries. Moreover, if the component of degree $2$ 
is non--zero, then the only possibility is that it maps 
$\fg_{-1} \times \fg_{-1}$ to $\fg_{0}$.  Thus it has its values in 
$ \fg_{-}^* \wedge \fg_{-}^* \otimes \fp$. It follows directly from 
the homogeneity reasons that components of degree $\geq 3$ have to have 
their values in this subbundle, too, and the geometry is torsion--free.

Finally, let us remind that vanishing of the harmonic curvature 
implies vanishing of the whole curvature, see Theorem \ref{harmonic}. 
This applies if $\ka_H$ coincides with $\ka^{(1)}$ which has to vanish 
for symmetric geometries, and they are locally isomorphic to 
homogeneous models. Now, it suffices to discuss components of 
harmonic curvature for concrete geometries, see Appendix.
\end{proof}

\subsection{Obstructions to flatness and involutive symmetries} 
\label{obstructions}
One can see from the discussion of the harmonic curvature that among 
all normal parabolic contact geometries, only 
\begin{itemize}
\item contact projective structures,
\item Lagrangean contact structures,
\item partially integrable almost CR--structures of hypersurface type
\end{itemize}
can carry a symmetry at a point with non--zero harmonic curvature. For each 
such symmetric geometry, there is exactly one obstruction to being 
locally isomorphic to the homogeneous model of the same type. 
For three--dimensional almost CR--structures and three--dimensional 
Lagrangean contact structures, there is the harmonic curvature 
$\ka^{(4)}$ valued in $\fg_{-1}^* \wedge \fg_{-2}^* \otimes \fg_1$. 
For the other ones, we have the harmonic curvature $\ka^{(2)}$ 
valued in $\fg_{-1}^* \wedge \fg_{-1}^* \otimes \fg_{0}$. 

Let us first focus on $\ka^{(2)}$.
Let $s_x$ be a symmetry at $x$ on a normal symmetric parabolic 
contact geometry and let $\va$ be as usual. For $u \in p^{-1}(x)$ 
and $X,Y \in \fg_{-1}$ we have
$$
\ka^{(2)}(\va(u))(X,Y)=\ka^{(2)}(u g_0 \exp Z)(X,Y)=
\exp(-Z) g_0^{-1} \cdot \ka^{(2)}(u)(X,Y),
$$
where $\cdot$ is the induced $\Ad$--action on 
$\fg_{-1}^* \wedge \fg_{-1}^* \otimes \fg_{0}$.  For the action 
of $g_0^{-1}$ we can write
$$
g_0^{-1} \cdot \ka^{(2)}(u)(X,Y)=
\Ad_{g_0^{-1}}(\ka^{(2)}(u)(\Ad_{g_0}X,\Ad_{g_0}Y))=
$$
$$
\Ad_{g_0^{-1}}(\ka^{(2)}(u)(-X,-Y))=\Ad_{g_0^{-1}}(\ka^{(2)}(u)(X,Y)).
$$ 
Because  $\fg_0$ is a subspace of 
$L(\fg_{-1},\fg_{-1}) \simeq \fg_{-1}^* \otimes \fg_{-1}$, 
the element $g_0^{-1}$ has to act trivially on $\fg_0$ and thus 
on $\ka^{(2)}(u)(X,Y)$ for each $X,Y$. Because also $\exp(-Z)$ acts 
trivially on $\ka^{(2)}(u)$, we get no additional restriction. In 
fact, $\ka^{(2)}$ is a tensor of type 
$\wedge^2T^{-1*}M \otimes T^{-1*}M \otimes T^{-1}M$ which is invariant 
with respect to the symmetry $s_x$. Let us denote this tensor by $W$.

We try to differentiate $W$ with respect to various Weyl connections. 
We focus on  connections corresponding to almost $s_x$--invariant 
Weyl structures, i.e. Weyl structures $\si$ satisfying 
$
\va(\si(u_0))=\si(\va_0(u_0)) \exp \U_1(u_0) \exp \U_2(u_0)
$
for suitable $\U_2$ and $\U_1$ such that $\U_1=0$ at $x$,
see \ref{action} for details.
\begin{lemma*}
On a symmetric normal parabolic contact geometry with a symmetry $s_x$ 
at $x$, let $\si$ be an arbitrary almost $s_x$--invariant Weyl structure 
and $\na^\si$ the corresponding Weyl connection. Then
\begin{enumerate}
\item[$(a)$] $\na^\si_{\xi} W=0$ holds at $x$ for each $\xi$ from the 
contact distribution,
\item[$(b)$] $\{\U_2, \xi\} \bullet W=0$ holds at $x$ for each  $\xi$ 
such that $T_xs_x.\xi(x)=\xi(x)$ and $\U_2$ is determined by 
$s_x$ at $x$.
\end{enumerate}
\end{lemma*}
\begin{proof}
In general,
$
\va(\si(u_0))=\si(\va_0(u_0)) \exp \U_1(u_0) \exp \U_2(u_0)
$
holds for each Weyl structure $\si$ and suitable $\U_1$ and $\U_2$, 
and this can be rewritten as
$$
\va^{-1}(\si(\va_0(u_0))) = \si(u_0)  \exp(-\U_1(u_0))  \exp(-\U_2(u_0)) 
$$
or simply $\va^*\si = \si \exp(-\U_1) \exp(-\U_2)$. 
For corresponding Weyl connections we then have 
$\na^{\va^*\si}=\na^{\si \exp(-\U_1) \exp(-\U_2)}$ and if we apply 
this on $W$, we get 
\begin{align} \label{nablaW}
\na^{\va^*\si}_\xi W=\na^{\si\exp(-\U_1)\exp(-\U_2)}_\xi W
\end{align}
for each vector field $\xi$. Because we suppose that $\si$ is almost 
$s_x$--invariant Weyl structure, then moreover $\U_1$ equals to zero in 
the fiber over $x$.

Let us discuss both sides of equation $(\ref{nablaW})$ at $x$ in detail. 
We start with the left hand side. At the point $x$ we have
$$
\na^{\va^*\si}_\xi W(\eta,\mu)(\nu)=(s_x^*\na^{\si})_\xi W(\eta,\mu)(\nu)= 
$$
$$
T_x s_x^{-1}. \na^{\si}_{T_x s_x.\xi} W(T_x s_x.\eta,T_x s_x.\mu)(T_xs_x.\nu)=
$$
$$
(-1)^4  \na^{\si}_{T_x s_x.\xi} W(\eta,\mu)(\nu)= \na^{\si}_{T_x s_x.\xi} W(\eta,\mu)(\nu)
$$
for each $\xi \in \frak X(M)$ and $\eta, \mu, \nu \in \Gamma(T^{-1}M)$ since  
$T_xs_x$ gives $-\id$ on $T^{-1}_xM$. 

Now we focus on the right hand side of (\ref{nablaW}). We use the 
general formula for the change of Weyl connections, see 
formula $(\ref{zmena-konexe})$ in \ref{WeylStructures}. Because our 
Weyl structure $\si$ satisfies $\U_1=0$ over $x$, the right hand side 
of $(\ref{nablaW})$ simplifies to 
$$
\na^{\si  \exp(-\U_1) \exp(-\U_2)}_\xi W = \na^\si_\xi W
+ \{\U_2, \xi_{-2}\} \bullet W
$$
in the fiber over $x$. If we put the above observations together, we see 
that the equation $(\ref{nablaW})$ can be rewritten as
\begin{align} \label{nablaW-upravena}
\na^{\si}_{T_x s_x.\xi} W = \na^\si_\xi W
+ \{\U_2, \xi_{-2}\} \bullet W
\end{align}
in the fiber over $x$. Let us discuss some concrete choices of the 
vector $\xi(x)$:
\smallskip\\
$(a)$ Suppose $\xi$ is contained in the contact 
distribution $T^{-1}M$. In particular, $\xi(x)=\xi_{-1}(x)$. 
Then $T_xs_x.\xi(x)=-\xi(x)$ and the equation 
$(\ref{nablaW-upravena})$ simplifies to
$$
\na^{\si}_{-\xi} W = \na^\si_\xi W
$$
at $x$. The algebraic bracket simply vanishes because $\xi_{-2}(x)=0$ 
in this case.
This implies $\na^\si_\xi W=0$ at $x$ for each $\xi$ from the 
contact distribution.
\smallskip\\
$(b)$ Let us now suppose that $T_xs_x.\xi(x)=\xi(x)$. Such vectors exist 
and for an almost $s_x$--invariant Weyl structure $\si$, these are 
exactly the  vectors satisfying $\xi=\xi_{-2}$ at $x$,  see \ref{action}. Then 
the equation $(\ref{nablaW-upravena})$ simplifies to
$$
\na^{\si}_{\xi} W = \na^\si_\xi W
+ \{\U_2, \xi\} \bullet W
$$
and we get the restriction
$\{\U_2, \xi \} \bullet W=0$
in the fiber over $x$. 
\end{proof}

The part $(a)$ is not surprising. Actually, $\na^\si W$ defines a tensor 
of type 
$T^{-1*}M \otimes \wedge^2 T^{-1*}M \otimes T^{-1*}M \otimes T^{-1}M$, 
i.e. of odd degree, which is invariant with respect to $s_x$. The 
consequences of the part $(b)$ we formulate in the following statement.

\begin{prop*}
On a symmetric normal parabolic contact geometry with a symmetry $s_x$ 
at $x$, suppose that $W$ is non--zero at $x$. Then each almost 
$s_x$--invariant Weyl structure has to be $s_x$--invariant. 
\end{prop*}  
\begin{proof}
With the above notation, we will discuss the formula $(b)$ from the Lemma 
for an almost $s_x$--invariant Weyl structure $\si$ and for some vector 
field $\xi$ such that $\xi(x)$ is non--zero and 
satisfies $T_xs_x.\xi(x)=\xi(x)$. Let us point out that such vectors exist 
and satisfy $\xi(x)=\xi_{-2}(x)$ in the isomorphism $TM \simeq \gr(TM)$ 
given by $\si$, see \ref{action}.
In some concrete frame $u =\si(u_0)$ from the fiber over $x$, we can 
write $\xi(x)=\lz u, X \pz$ for suitable non--zero $X \in \fg_{-2}$. 
Similarly, $\U_2(x)=\lz u, Z \pz$ for suitable $Z \in \fg_2$ and the 
algebraic bracket $\{\U_2, \xi \}$ corresponds to $\lz u, [Z,X] \pz$ at $x$.
Moreover, if $Z \neq 0$, we can choose $X$ such that $[Z,X]$ is exactly 
the grading element $E$, see \ref{contact}. Then $\{\U_2, \xi \}$ 
corresponds to the grading section $E(x)$, see \ref{WeylStructures}. 
In particular, it acts by the algebraic action $\bullet$ on $W$ by 
its homogeneity. Because $W$ has homogeneity two, we get 
$$\{\U_2, \xi \}\bullet W =2W$$
in the fiber over $x$ and we have a restriction of the form $2W=0$ at 
$x$. This is a contradiction with the assumption that $W$ is non--zero 
at $x$. Thus  the only possibility is that $Z=0$ and thus $\U_2$ has to 
vanish at $x$.  But this means that almost $s_x$--invariant Weyl 
structure is $s_x$--invariant, see \ref{action}.
\end{proof}

Before we proceed further, let us return to $\ka^{(4)}$ valued in 
$\fg_{-1}^* \wedge \fg_{-2}^* \otimes \fg_1$. The discussion of $\ka^{(4)}$ 
is parallel to the discussion of $\ka^{(2)}$ and we summarize it very 
briefly. If $\ka^{(4)}$ is non--zero, then it defines a tensor of even 
degree which has to be invariant with respect to the symmetry $s_x$.
Really, for each $X \in \fg_{-1}$ and $V \in \fg_{-2}$ we have
$$
\ka^{(4)}(\va(u))(X,V)=\ka^{(4)}(u g_0 \exp Z)(X,V)=
\exp(-Z) g_0^{-1} \cdot \ka^{(4)}(u)(X,V)=
$$
$$
\Ad_{g_0^{-1}}(\ka^{(4)}(u)(\Ad_{g_0}X,\Ad_{g_0}V))= -\ka^{(4)}(u)(-X,V)=\ka^{(4)}(u)(X,V).
$$
Let us denote this tensor by $Y$. Again, we can differentiate $Y$ with 
respect to a Weyl connections corresponding to almost invariant Weyl structures.
For an almost $s_x$--invariant Weyl structure $\si$ we have the equation
$$
\na^{\va^*\si}_\xi Y=\na^{\si \exp(-\U_1) \exp(-\U_2)}_\xi Y
$$
for each vector field $\xi$ and suitable $\U_1$ and $\U_2$ corresponding 
to $\si$, where $\U_1$ vanishes at $x$. The left hand side 
can be rewritten as
$$
\na^{\va^*\si}_\xi Y(\eta,\mu)=(s_x^*\na^{\si})_\xi Y(\eta,\mu)= 
T_x s_x^{-1}. \na^{\si}_{T_x s_x.\xi} Y(T_x s_x.\eta,T_x s_x.\mu)
$$
$$
=(-1)^2 \na^{\si}_{T_x s_x.\xi} Y(\eta,\mu)= \na^{\si}_{T_x s_x.\xi} Y(\eta,\mu)
$$
for each $\xi \in \frak X(M)$, $\eta \in \Gamma(T^{-1}M)$ and 
$\mu \in \frak X(M)$ such that $\mu=\mu_{-2}$ via the isomorphism 
given by $\si$ at $x$.
Really, $T_xs_x$ gives $-\id$ on $T^{-1}_xM$ and $T_xs_x.\mu(x)=\mu(x)$.
Thus we get the restriction of the form 
$$
\na^\si_{T_xs_x.\xi} Y = \na^\si_{\xi} Y + \{\U_2, \xi_{-2}\} \bullet Y
$$
and we have $\na^\si_{\xi_{-1}}Y=0$ and $\{\U_2, \xi_{-2}\} \bullet Y=0$ 
in the fiber over $x$. Because $Y$ is of homogeneity four, the same 
arguments as in the proof of the above Proposition shows that $\U_2$ 
vanishes at $x$ and then each $s_x$--invariant Weyl structure has to be 
$s_x$--invariant. All these observations together with the last 
Proposition give us the following statement.

\begin{thm*}
On a symmetric normal parabolic contact geometry with a symmetry $s_x$ at 
$x$, suppose that its harmonic curvature is non--zero at $x$. Then $s_x$ is involutive. 
\end{thm*}
\begin{proof}
In such case, almost $s_x$--invariant Weyl structures have to be 
$s_x$--invariant and the rest follows from  \ref{action}. 
\end{proof}
\begin{cor*}
On a symmetric normal parabolic contact geometry with a symmetry $s_x$ 
at $x$, suppose that its harmonic curvature is non--zero at $x$. Then there 
are admissible affine connections which are invariant with respect to 
the symmetry $s_x$: We take Weyl connections corresponding 
to $s_x$--invariant Weyl structures.
\end{cor*}

\section{Uniqueness of symmetries}
We discuss here the question how many different symmetries can exist at 
a point with non--zero curvature. We first give one general restriction 
and then some consequences for concrete geometries. 

\subsection{Algebraic restriction} \label{algebraic}
Let $s_x$ and $\bar s_x$ are two different symmetries at $x$ on a 
symmetric normal parabolic contact geometry with non--zero harmonic curvature 
at $x$ and denote by $\va$ and $\bar \va$ corresponding automorphisms 
of the parabolic geometry. Clearly, $s_x \neq \bar s_x$ if and only if 
$\va \neq \bar \va$. Then symmetries $s_x$ and $\bar s_x$ are involutive 
and there exist $s_x$--invariant and $\bar s_x$--invariant Weyl 
structures, see \ref{obstructions} and \ref{action}.

\begin{lemma*}
For each two different involutive symmetries $s_x$ and $\bar s_x$ at $x$ on a 
symmetric  parabolic contact geometry,  $s_x$--invariant and $\bar s_x$--invariant Weyl structures 
form two disjoint families of Weyl structures. 
\end{lemma*}
\begin{proof}
Suppose there is a Weyl structure $\si$ which  is $s_x$--invariant and 
$\bar s_x$--invariant at $x$, i.e. $\va(\si(u_0)) = \si(\va_0(u_0))$ 
and simultaneously $\bar \va (\si(u_0))=\si(\bar \va_0(u_0))$ in the fiber
 over $x$. Then, the corresponding Weyl connection $\na^\si$ is invariant 
with respect to both symmetries $s_x$ and $\bar s_x$. But similarly as in the 
last part of the proof of Proposition \ref{action}, the connection 
$\na^\si$ determines uniquely the symmetry via behavior of its 
geodesics at $x$. Consequently, $s_x=\bar s_x$ on a neighborhood of $x$. 
\end{proof}
Let $\si$ be $s_x$--invariant Weyl structure and let $\bar \si$ be 
$\bar s_x$--invariant Weyl structure. Then 
$\bar \si=\si \exp \U_1 \exp \U_2$ holds for suitable 
$\U_1:\ba_0 \rightarrow \fg_1$ and $\U_2:\ba_0 \rightarrow \fg_2$. 
The last Lemma says that $\U_1$ has to be non--zero at $x$. 
\begin{prop*} 
Suppose there are two different involutive symmetries at $x$ on a symmetric 
normal parabolic contact geometry 
and let $\si$ and $\bar \si$ 
are corresponding invariant Weyl structures. For all $\xi$ from the 
contact distribution, the bracket $\{\U_1, \xi\}$ acts trivially by 
the algebraic action on $W$ or $Y$, respectively, at $x$.
\end{prop*}
\begin{proof}
Let us start with $W$. Let $\xi$ be an arbitrary vector field from the 
contact distribution, thus $\xi=\xi_{-1}$ for each Weyl structure. 
The Lemma \ref{obstructions} gives $\na^\si_{\xi_{-1}} W=0$ 
 and $\na^{\bar \si}_{\xi_{-1}}W=0$ at $x$. Simultaneously, we have 
$\bar \si=\si  \exp \U_1 \exp \U_2$ and the formula 
$(\ref{zmena-konexe})$ from \ref{WeylStructures} gives
$$
\na^{\bar \si}_{\xi_{-1}}W=\na^{\si}_{\xi_{-1}}W + \{ \xi_{-1}, \U_1\} \bullet W
$$
at $x$, since $\xi_{-2}(x)=0$. Because both covariant derivatives vanish 
at $x$, we get the restriction of the form $\{ \xi_{-1}, \U_1\} \bullet W=0$
at $x$ for each $\xi$ from the contact distribution.
One can see from \ref{obstructions} that the same line of arguments works 
for $Y$ and we get the restriction of the form 
$\{ \xi_{-1}, \U_1\} \bullet Y=0$ at $x$ for each $\xi$ from the 
contact distribution.
\end{proof}
\begin{rem*}
Let us again point out that the existence of a non--involutive symmetry at $x$ causes vanishing of the harmonic curvature at $x$, see \ref{obstructions}. 
\end{rem*}
\subsection{Examples}
Let us now discuss the above restrictions for concrete types of geometries. 
The key point is to find sufficiently nice $\xi$ such that the action of 
the above algebraic bracket is easily understandable. 

\subsubsection*{Lagrangean contact structures} 
Let us first point out that we use here the notation from \ref{examples}. 
The decomposition of the contact distribution into two isotropic 
subbundles $T^{-1}M=L \oplus R$ can be interpreted as a product structure 
on $T^{-1}M$, which an operator $J:T^{-1}M \rightarrow T^{-1}M$ satisfying $J^2=\id$. The subbundles $L$ and $R$ are simply eigenspaces of $J$. 
The Levi bracket $\Le: T^{-1}M \times T^{-1}M \rightarrow TM/T^{-1}M$ 
is non--degenerate antisymmetric bilinear map, and then, $\Le(- , J- )$ 
is non--degenerate symmetric map which defines a conformal class 
of pseudometrics on $T^{-1}M$ of signature $(n,n)$. We denote the class by $g$. 
Each pseudometric is then given by the choice the identification 
$TM/T^{-1}M \simeq \R$. In particular,
the question whether $g(\xi,\eta)$ equals to zero for some 
$\xi,\eta \in T^{-1}M$ makes sense, because the answer does not depend on 
the choice of the metric from the class.
  
\begin{prop*}
Suppose there are two different involutive symmetries at $x$ on a symmetric 
normal Lagrangean contact structure and denote by $\si$ and 
$\si \exp \U_1 \exp \U_2$  corresponding invariant Weyl structures. 
Identify $\U_1$ with its image in $T^{-1}M=L\oplus R$ via an isomorphism 
given by a metric from $g$ and denote by $\U_1^L$ and $\U_1^R$ 
corresponding components in $L$ and $R$. If $g(\U_1^R,\U_1^L) \neq 0$ at 
$x$, then the harmonic curvature vanishes at $x$. 
\end{prop*}
\begin{proof}
We discuss the restriction from the Proposition \ref{algebraic} for 
Lagrangean contact structures in detail.
Let us write $\U_1(x)=\lz u,Z\pz$ for suitable
$$
Z= \Big( \begin{smallmatrix}
0 & S & 0 \\ 0 & 0 & T \\ 0 & 0 & 0
\end{smallmatrix} \Big) \in \fg_{1}
,$$
which has to be non--zero, see Lemma \ref{algebraic}.
Choose $\xi_{-1} \in \Gamma(T^{-1}M)$ such that 
$\xi_{-1}(x)=\lz u, X \pz$ for  $X$ of the form
$$
X=\Big( \begin{smallmatrix}
0 & 0 & 0 \\ T & 0 & 0 \\ 0 & S & 0
\end{smallmatrix} \Big) \in \fg_{-1}
.$$ 
The bracket $\{ \xi_{-1},\U_1\} $ then corresponds to 
$\lz u, [X,Z] \pz$ at $x$, where 
$$
[X,Z]=  \Big( \begin{smallmatrix}
-ST & 0 & 0 \\ 0 & 0 & 0 \\ 0 & 0 & ST
\end{smallmatrix} \Big) \in \fg_{0}.
$$
It is easy to verify that with this choice, $[X,Z]$ is simply a 
grading element multiplied by a non--zero number $-ST$. Then the 
bracket $\{ \xi_{-1}, \U_1 \}$ is a non--zero multiple of the grading 
section $E(x)$. Via the identification given by the metric from $g$, 
the components $S$ and $T$ correspond to components $\U_1^R$ and $\U_1^L$ 
of $\U_1$ in subbundles $R$ and $L$ at $x$ and the fact that $ST\neq 0$ 
means that $g(\U_1^R,\U_1^L) \neq 0$ at $x$.
Because the grading section acts on $W$ by its homogeneity, 
$\{ \xi_{-1}, \U_1 \}$ acts trivially on $W$ if and only $W$ vanishes 
at $x$. Clearly, the same arguments work for $Y$.
\end{proof}

\subsubsection*{Partially integrable almost CR--structures} 
Let us first point out that we use here the notation from 
\ref{examples}. Moreover, suppose that the geometry is oriented and then, 
we can speak about signature of the structure. Using the complex 
structure $J$ given on $T^{-1}M$, we can define a non--degenerate 
symmetric mapping $\Le(-,J-)$, which defines a conformal class 
of pseudometrics on $T^{-1}M$. The signature is given by the signature 
of the structure. Let us denote the class by $g$. Each 
pseudometric from the class is given by the choice 
of the identification $TM/T^{-1}M \simeq \R$. In particular, the 
question whether $g(\xi,\xi) \neq 0$ for $\xi \in T^{-1}M$ makes 
sense, because the answer does not depend on the choice of the 
pseudometric from the class.

\begin{prop*}
Suppose there are two different involutive symmetries at $x$ on a 
symmetric normal partially integrable almost CR--structure  
and denote by  $\si$ and 
$\si \exp \U_1 \exp \U_2$  corresponding invariant Weyl 
structures.  Identify $\U_1$ with its image in $T^{-1}M$ via an 
isomorphism given by a metric from $g$. If $g(\U_1,\U_1) \neq 0$ at $x$, 
i.e. if the length of $\U_1$ is non--zero at $x$, then the harmonic curvature 
vanishes at $x$.
\end{prop*}
\begin{proof}
We discuss the restriction from the Proposition \ref{algebraic} 
for CR--structures in detail. Let us write $\U_1(x)=\lz u,Z\pz$ for suitable
$$
Z= \Big( \begin{smallmatrix}
0 & 0 & 0 \\ IS^* & 0 & 0 \\ 0 & -S & 0
\end{smallmatrix} \Big)
,$$  
which has to be non--zero, see Lemma \ref{algebraic}. Choose 
$\xi_{-1} \in T^{-1}M$ such that $\xi_{-1}(x)=\lz u, X \pz$ for  $X$ of the 
form
$$
X=
\Big( \begin{smallmatrix}
0 & 0 & 0 \\ S & 0 & 0 \\ 0 & -S^*I & 0
\end{smallmatrix} \Big)
\in \fg_{1}
.$$ 
The bracket $\{ \xi_{-1},\U_1\} $ then corresponds to $\lz u, [X,Z] \pz$ 
at $x$, where 
$$
\Big( \begin{smallmatrix}
-SIS^* & 0 & 0 \\ 0 & 0 & 0 \\ 0 & 0 & SIS^*
\end{smallmatrix} \Big) \in \fg_{0}.
$$
It is easy to verify that with this choice, $[X,Z]$ is simply a 
grading element multiplied by a non--zero number $-SIS^*$. Then the 
bracket $\{ \xi_{-1}, \U_1 \}$ is a non--zero multiple of the grading 
section $E(x)$. 
Using the identification $T^{-1}M \simeq T^{-1*}M$ given by a metric 
from $g$, $-SIS^*$ corresponds to $g(\U_1,\U_1)$ and $SIS^*\neq 0$ 
means that $g(\U_1,\U_1) \neq 0$. Because the grading section acts on $W$ 
by its homogeneity, $\{ \xi_{-1}, \U_1 \}$ acts trivially on $W$ if and 
only $W$ vanishes at $x$. Clearly, the same arguments work for $Y$.
\end{proof}

\begin{cor*}
Suppose there are two different involutive symmetries at $x$ on a 
symmetric normal strictly pseudoconvex partially integrable 
almost CR--structure. Then the harmonic curvature vanishes at $x$. 
\end{cor*}

\section*{Appendix: Contact gradings and corresponding geometries} \label{tabulka}
Let us sketch here briefly a classification of contact gradings of 
real semisimple Lie algebras. There is the well know classification of 
all (complex) semisimple Lie algebras in the language Dynkin diagrams 
and description of all their real forms in the language of Satake 
diagrams, see \cite{parabook,Y}. 
It can be proved that if a Lie algebra admits a contact grading, then it 
has to be simple. It turns out that except 
$\frak{sl}(2,\R)$, $\frak{sl}(n,\mathbb{H})$, 
$\frak{so}(n-1, 1)$, $\frak{sp}(p, q)$ and 
some real forms of $E_6$ and $F_4$, any non--compact non--complex real 
simple Lie algebra admits a unique real contact grading, see \cite{parabook}. 

Let us start with real classical Lie algebras, i.e. real forms of Lie 
algebras of type $A_\ell$, $B_\ell$, $C_\ell$ and $D_\ell$. In the first 
row of the following table, we indicate a real simple Lie algebra which 
admits a contact grading. In the second row we specify the geometry, 
which corresponds to the unique contact grading and in the last row we 
write its components of harmonic curvature.
\vspace{2mm}

\begin{center}
\begin{tabular}{|p{2,5cm}|p{5cm}|p{3cm}|}
\hline
real simple $\fg$ & contact geometry &  components of $\ka_H$ 
\\
\hline \hline
 $\frak{sl}(3,\R)$ 
& Lagrangean contact structures of dimension  $3$ & 
$\fg_{-2} \times \fg_{-1} \rightarrow \fg_{1}$ 
$\fg_{-2} \times \fg_{-1} \rightarrow \fg_{1}$ 
\\ \hline
 $\frak{sl}(n+2,\R)$ for $n\geq 2$ & Lagrangean contact structures of dimension $2n+1$ & 
$\fg_{-1} \times \fg_{-1} \rightarrow \fg_{-1}$  
$\fg_{-1} \times \fg_{-1} \rightarrow \fg_{-1}$ 
$\fg_{-1} \times \fg_{-1} \rightarrow \fg_{0}$  
\\ \hline
$\frak{su}(2,1)$ and $\frak{su}(1,2)$ & partially integrable almost CR structures of dimension $3$ & 
$\fg_{-2} \times \fg_{-1} \rightarrow \fg_{1}$ 
\\ \hline
$\frak{su}(p+1,q+1)$ for $p+q \geq 2$ & partially integrable almost CR structures of dimension $2p+2q+1$ & 
$\fg_{-1} \times \fg_{-1} \rightarrow \fg_{-1}$   
$\fg_{-1} \times \fg_{-1} \rightarrow \fg_{0}$  
\\ \hline
$\frak{so}(p+2,q+2)$  $p+q \neq 4$ & Lie contact structures of dimension $2p+2q+1$ & $\fg_{-1} \times \fg_{-1} \rightarrow \fg_{-1}$  $\fg_{-1} \times \fg_{-1} \rightarrow \fg_{-1}$
\\ \hline
$\frak{so}(p+2,q+2)$ for $p+q =4$ & Lie contact structures of dimension $9$ & $\fg_{-1} \times \fg_{-1} \rightarrow \fg_{-1}$  $\fg_{-1} \times \fg_{-1} \rightarrow \fg_{-1}$ $\fg_{-1} \times \fg_{-1} \rightarrow \fg_{-1}$
 \\ \hline
$\frak{sp}(n+2)$ for $n \geq 1$ & contact projective structures & $\fg_{-1} \times \fg_{-1} \rightarrow \fg_{0}$ \\ \hline
\end{tabular}
\end{center}
\vspace{2mm}

Let us also give a brief overview of contact gradings corresponding 
to exotic Lie algebras. 
For types $G_2$ and $F_4$, there is exactly one real algebra admitting 
contact grading, the split real form. For $E_6$, there are three real 
forms which admit a contact grading, the split form and two 
$\frak{su}$--algebras.  For $E_7$, there are three different real forms 
and for $E_8$, there are two different real forms admitting a contact 
grading. The description of corresponding geometries can be found 
in \cite{parabook}. All these geometries have harmonic curvatures only of 
type $\ka^{(1)}$ valued in $ \fg_{-1}^* \wedge \fg_{-1}^* \otimes \fg_{-1}$.

\end{document}